\documentclass[12pt]{article}
\usepackage[final]{epsfig}
\usepackage{graphics}
\usepackage{amsmath}
\usepackage{amsfonts}
\usepackage{latexsym}
\usepackage{amssymb}
\usepackage{graphicx}
\usepackage{url}
\usepackage{epstopdf}

\newtheorem{lemma}{Lemma}[section]

\newtheorem{remark}[lemma]{Remark}
\newtheorem{example}[lemma]{Example}
\newtheorem{theorem}{Theorem}

\begin{document}
\newcommand{\eps}{{\varepsilon}}
\newcommand{\proofend}{$\Box$\bigskip}
\newcommand{\C}{{\mathbb C}}
\newcommand{\Q}{{\mathbb Q}}
\newcommand{\R}{{\mathbb R}}
\newcommand{\Z}{{\mathbb Z}}
\newcommand{\RP}{{\mathbf {RP}}}
\newcommand{\CP}{{\mathbf {CP}}}
\def\proof{\paragraph{Proof.}}

\title{On  Lagrangian tangent sweeps and Lagrangian outer billiards}

\author{Dmitry Fuchs\footnote{
Department of Mathematics, University of California, Davis, CA 95616;
 fuchs@math.ucdavis.edu} \and Serge Tabachnikov\footnote{
Department of Mathematics,
Pennsylvania State University,
University Park, PA 16802;
tabachni@math.psu.edu}
}

\date{}
\maketitle

\section{Introduction: tangent sweeps of plane curves and planar outer billiards}
\begin{figure}[hbtp]
\centering
\includegraphics[height=1.5in]{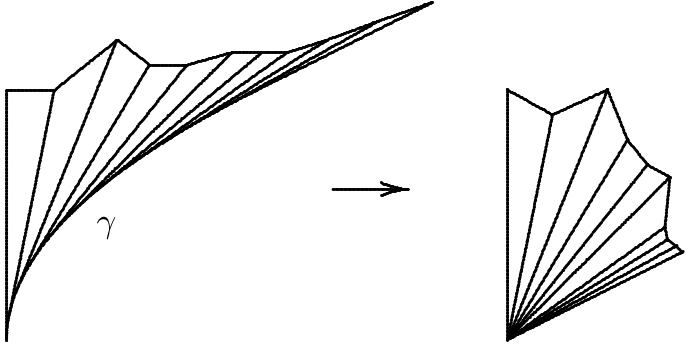}
\caption{An area preserving  map}
\label{segments}
\end{figure}

The exterior of a convex oriented smooth plane curve $\gamma$ is foliated by the tangent rays to the curve. Fix an origin and parallel translate each tangent ray so that its endpoint moves to the origin. This defines a map from the exterior of the curve to the plane. It is a calculus exercise to prove that this map is area preserving. See Figure \ref{segments}.

\begin{figure}[hbtp]
\centering
\includegraphics[height=1.45in]{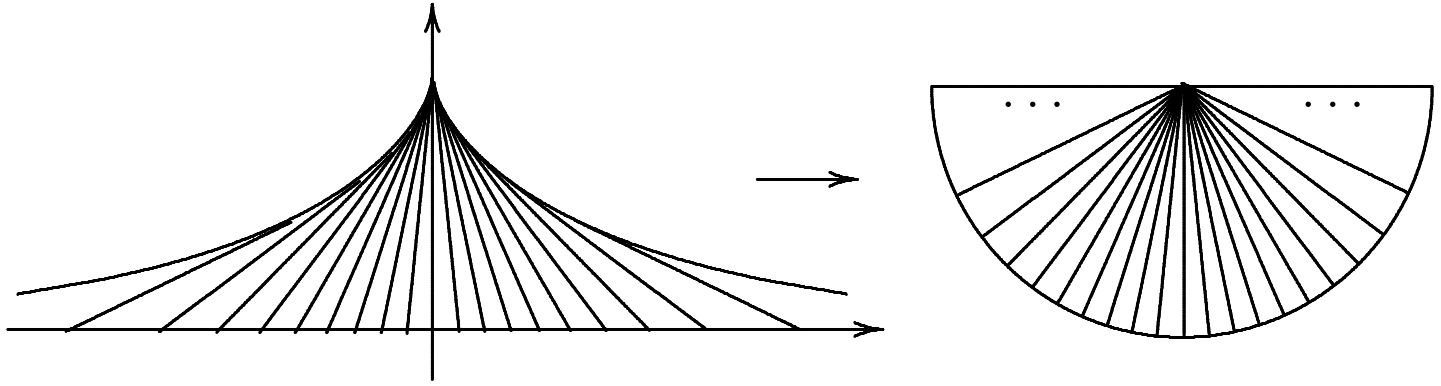}
\caption{The tractrix}
\label{tractrix}
\end{figure}

In particular, this makes it possible to calculate the areas swept by the tangent segments to a curve. For example, the classical tractrix is defined by the property that the  locus of the endpoints of its unit tangent segments is a straight line. The tangent segment to the tractrix makes half a turn, and the area under the tractrix equals $\pi/2$, see Figure \ref{tractrix}. 

In the terminology of \cite{AM}, ``the tangent sweep of a curve and its tangent cluster have the same area".  We refer to \cite{AM} for numerous applications of this observation, which is referred to therein as ``Mamikon's Theorem".

Changing the orientation of the curve $\gamma$ yields another area preserving map. Combined with the fact that the reflection in the origin is area preserving, this
implies that the following {\it outer billiard transformation}\footnote{The notion of outer billiards appears in the mathematical literature under different names; in particular, in all works mentioned in our bibliography outer billiards are called {\it dual billiards}.} is area preserving. 

\begin{figure}[hbtp]
\centering
\includegraphics[height=1.7in]{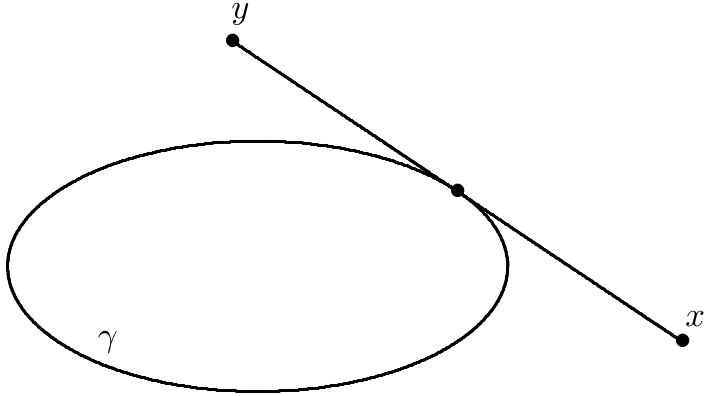}
\caption{Outer billiard map about $\gamma: x \mapsto y$}
\label{bill}
\end{figure}

Two points are in the outer billiard correspondence with respect to a curve $\gamma$ if they lie on the same tangent line to $\gamma$ at equal distances from the tangency point. If $\gamma$ is an oriented smooth closed strictly convex curve, this correspondence is an area preserving diffeomorphism of the exterior of $\gamma$, see Figure \ref{bill}.

Starting with J. Moser \cite{Mo}, outer billiards have been an active area of research, see \cite{DT,Ta2} for surveys. 

The plane is a 2-dimensional linear symplectic space, and a curve is its Lagrangian submanifold. In what follows, we extend the `tangent sweep' area preserving property and the outer billiard correspondence to higher-dimensional symplectic spaces and their Lagrangian submanifolds.

\section{Lagrangian tangent sweep property}

Consider linear symplectic space $(\R^{2n},\omega)$ and a Lagrangian submanifold $L^n$ in it. One can identify the symplectic space with the cotangent bundle $T^* \R^n$ with the usual $(q,p)$-coordinates, so that $\omega=dp\wedge dq$, and so that $L\subset T^* \R^n$ is locally the graph of the differential of a (generating) function $F(q_1,\dots,q_n)$:
\begin{equation}
\label{gener}
p_i=F_{q_i},\ i=1,\ldots,n.
\end{equation}
Here and elsewhere we use the shorthand notations $F_{q_i}=\partial F/ \partial {q_i}$,  $dp\wedge dq=\sum dp_i \wedge dq_i$, etc.

In this section, we are interested in local questions, and we assume that $L$ is as in (\ref{gener}).
We refer to \cite{AG} for an introduction to symplectic geometry.

Adopting the terminology of \cite{AM}, the tangent sweep of the Lagrangian manifold $L$ is the union of its tangent spaces (considered as $n$-dimensional affine subspaces of $\R^{2n}$), and the  tangent cluster of $L$ is the result of parallel translating these spaces so that the foot point of each tangent space becomes the origin. This defines a multivalued map $\Phi$ from the tangent sweep to the  tangent cluster.

More precisely, the inclusion $L \subset \R^{2n}$ extends to the map $\varphi:TL\to \R^{2n}$. One also has the map $\psi: TL\to \R^{2n}$, parallel translating affine tangent spaces to the origin. Assume that $\varphi$ is locally a diffeomorphism. Then $\Phi=\psi\circ\varphi^{-1}$.

\begin{theorem}
\label{sympl}
The map $\Phi$ is a local symplectomorphism.
\end{theorem}

\proof
The tangent space at point $(q,p)\in L$ is given by
\begin{equation}
\label{tangs}
(q_1+t_1,\ldots,q_n+t_n,F_{q_1}+\sum t_i F_{q_1 q_i},\ldots,F_{q_n}+\sum t_i F_{q_n q_i}),
\end{equation}
where $(t_1,\ldots,t_n)$ are coordinates in the tangent space. (Here and elsewhere summation is over repeated indices).

Let $(X,Y)$ and $(x,y)$ be Darboux coordinates in the source and target spaces.
Then map $\Phi: (X,Y) \mapsto (x,y)$ is given by the formulas
\begin{equation}
\label{formulas}
x_i=t_i, y_i=\sum t_k F_{q_i q_k}, X_i=q_i+t_i=q_i+x_i, Y_i=F_{q_i}+\sum t_k F_{q_i q_k}=F_{q_i}+y_i.
\end{equation}
We want to check that $\sum dX_i\wedge dY_i = \sum dx_i\wedge dy_i$. Indeed,
$$
dX_i\wedge dY_i - \sum dx_i\wedge dy_i= \sum dq_i\wedge dF_{q_i} + \sum dq_i\wedge dy_i + \sum dx_i\wedge dF_{q_i}$$
$$=\sum dq_i\wedge dF_{q_i} + \sum t_k dq_i\wedge dF_{q_i q_k} + \sum F_{q_i q_k} (dq_i\wedge dt_k + dx_i\wedge dq_k).$$
The first sum vanishes because it equals $d^2(F)$, the second sum vanishes because each summand for fixed $k$ equals $t_k d^2(F_{q_k})$, and the third sum vanishes because $x_i=t_i$. 
\proofend

\section{Local multiplicity of covering by tangent spaces}
\label{localmult}

Denote by $\Delta$ the set of critical values of the map $\varphi$. To describe $\Delta$,
 we compute the Jacobian  of  $\varphi$. Let $A(q,t)$ be the $n\times n$ symmetric matrix given by
$$
a_{ij}=\sum t_k F_{q_i q_j q_k}.
$$

\begin{lemma}
\label{Jac}
Up to a sign, the Jacobian  of $\varphi$ equals $\det A$.
\end{lemma}

\proof
According to Theorem \ref{sympl}, we need to compute the volume form $(\sum dx_i\wedge dy_i)^n$. Using (\ref{formulas}), and ignoring the signs, this volume form equals
$$
dt_1\wedge\ldots\wedge dt_n\wedge dy_1\wedge\ldots\wedge dy_n =
\det A\ dt_1\wedge\ldots\wedge dt_n\wedge dq_1\wedge\ldots\wedge dq_n,
$$
and we are done.
\proofend

It follows that $\Delta$ consists of $L$ and, in every tangent space $T_{(q,p)}L$ to $L$, of a conical hypersurface of degree $n$ given by the homogeneous equation $\det A =0$. If $n\ge2$, then this hypersurface may be empty, as the next example shows, and it may be the whole tangent space (which does not occur in generic points of a generic $L$). The case $n=1$ is simple: $\Delta$ is always the union of $L$ and all tangent lines at the inflection points.

\begin{example} \label{Jacob}
{\rm Let
$$
F(q_1,q_2)=aq_1^2q_2+bq_1q_2^2,\ \ a,b\neq 0.
$$
Then $\det A\neq 0$ for all $(t_1,t_2) \neq (0,0)$.
}
\end{example}

\begin{figure}[hbtp]
\centering
\includegraphics[height=2.5in]{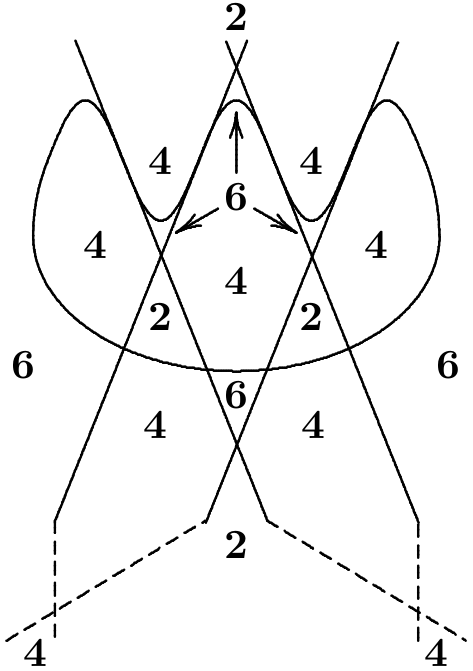}
\caption{The number of tangent spaces in the complement of $\Delta$.}
\label{curve}
\end{figure}

We are interested in the number of tangent spaces to $L$ that pass through a given point. This is a locally constant function in the complement to $\Delta$, see Figure \ref{curve}. In particular, if $\Delta = L$ and $n\ge 2$, as in the example above, then this function has the same value everywhere in the complement to $L$. In general, however, this function is unbounded (which is also clear from Figure \ref{curve}).

Consider a local version of this question.
Let $L\subset \R^{2n}$ denote a germ of a Lagrangian manifold defined by a generating function (\ref{gener}). 
We assume that the generating function $F(q)$ is real analytic and generic. One of the general position assumptions is that  the cubic part of $F$ contains  the terms $q_1^3,\ldots, q_n^3$ with non-zero coefficients. Other general position condition will be specified in the proof. 

We are interested in the multiplicity of the covering of $\R^{2n}$ by the tangent spaces of $L$.

\begin{theorem}
\label{multi}
The multiplicity in question is not greater than $2^n$, and this estimate is sharp.
\end{theorem}

\begin{figure}[hbtp]
\centering
\includegraphics[height=2.7in]{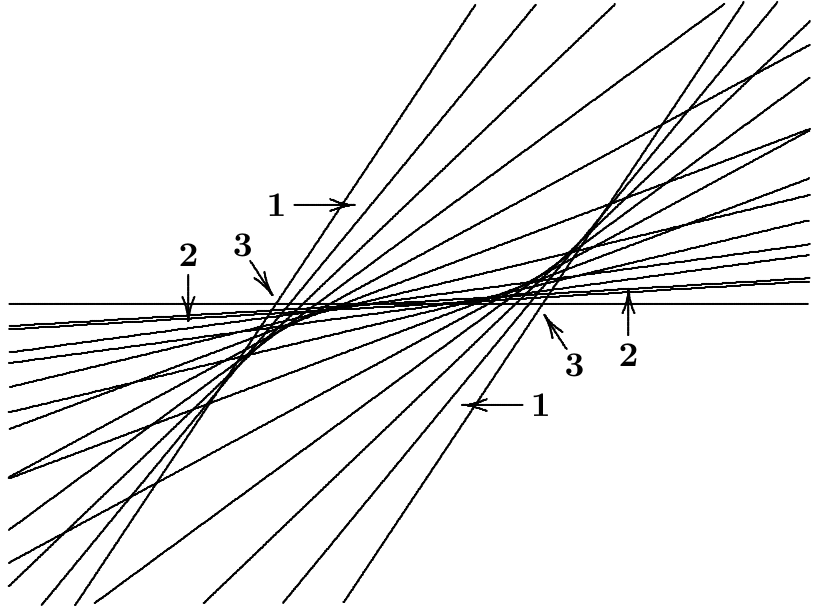}
\caption{The local multiplicity in a neighborhood of an inflection point.}
\label{inflection}
\end{figure}

\begin{remark}
{\rm As we have seen before, globally no estimate exists. Notice that even locally the estimate in Theorem does not hold, if the germ is not generic. For example, as Figure \ref{inflection} shows, the multiplicity may be greater than 2 in an arbitrarily small neighborhood of a curve with a (generic) inflection point in the plane.}
\end{remark}

\noindent{\bf Proof of Theorem} Let $(x,y)\in \R^{2n}$ be a test point. It follows from formula (\ref{tangs}) that 
$(x,y)\in T_{(q,p)}L$ if and only if there exist $(q,t)$ such that 
$$
x_i=q_i+t_i,\ y_i=F_{q_i}+\sum t_j F_{q_i q_j},\ i=1,\dots,n.
$$
One expresses $t$ from the first equation and substitutes to the second one:
$$
\sum x_j F_{q_i q_j} - \sum q_j F_{q_i q_j} + F_{q_i} =y_i, \ i=1,\dots,n,
$$
or 
\begin{equation}
\label{crit}
\nabla (x \cdot F_q - q\cdot F_q + 2F) =y.
\end{equation}

Without loss of generality, assume that $L$ is a germ at the origin and that the tangent space to $L$ at the origin is the $q$-space. Let $G$ be the cubic part of $F$ and $H$ be the terms of higher order: $F=G+H$.
Then the function on the left hand side of (\ref{crit}) is
$$
\Phi_x (q):= x \cdot G_q + (x\cdot H_q - G) + (2H - q\cdot H_q),
$$
(we used the Euler formula $q\cdot G_q = 3 G$). The three terms are quadratic, cubic, and higher order terms in $q$.

Consider the complexification of (\ref{crit}). The left hand side of this equation defines a germ of a holomorphic map $\nabla \Phi_x:  (\C^n,0) \to (\C^n,0)$. We are interested in the number of preimages in a neighborhood of the origin of a generic point $y$ near the origin. 

This number is given by the dimension of the local algebra of the map $\nabla \Phi_x$, see Section 4.3 of \cite{AGV}. Recall the definition.

Let ${\cal O}$ be the ring of germs of holomorphic functions of $n$ variables, and let
$\varphi=(\varphi_1,\dots,\varphi_n): (\C^n,0) \to (\C^n,0)$ be a germ of a holomorphic map. The local algebra of $\varphi$ is defined as the quotient algebra $Q_\varphi = {\cal O}/(\varphi_1,\dots,\varphi_n). $

In our situation, the components of the map are the first partial derivatives of the function $\Phi_x$. In this case, the dimension of the local algebra is called the Milnor number of the function, denoted by $\mu(\Phi_x)$.

To compute the Milnor number, we use Kouchnirenko's theorem \cite{Ko}.
According to this theorem, the Milnor number $\mu(f)$ of an isolated singularity of a generic germ of an analytic function $f$ of $n$ variables can be computed in terms of its Newton polyhedron (the genericity is understood as the complement to an algebraic subvariety in the space of coefficients). 

Let $\Gamma$ be the boundary of the Newton polyhedron of $f$, and $S$ be the union of the segments 
in $\R^n$ joining the origin to points of $\Gamma$. The Newton number is 
$$
\nu(f):=n!V_n-(n-1)!V_{n-1}+\ldots +(-1)^{n-1}1!V_1+(-1)^n,
$$
where $V_n$ is the $n$-dimensional volume of $S$ and, for $1\le k\le n-1$, $V_k$ is the sum of the $k$-dimensional volumes of the intersection of $S$ with the coordinate planes of dimension $k$. The statement is that $\mu(f)=\nu(f)$.

In our case, $f= \Phi_x$. If $x=0$ then $\Phi_x$ has no quadratic terms, and the Newton polyhedron $S$  is an $n$-dimensional simplex whose vertices are the origin and the points on the coordinate axes at distance 3 from the origin. The number of the simplices that are the intersections of $S$ with the coordinate planes of dimension $k$ is ${n\choose k}$, and the $k$-dimensional volume of each such intersection  equals $3^k/k!$ Accordingly, 
$$
\nu(G) = \sum_{k=0}^n (-1)^{n-k} {n\choose k} 3^k = (3-1)^n=2^n.
$$

Likewise, if the vector $x$ has $m>0$ non-zero components then the generic Newton polyhedron is an $n$-dimensional simplex whose vertices are the origin and the points on the coordinate axes, $m$ of them at distance 2, and $n-m$ at distance 3 from the origin. Its Newton number is  equal to $(2-1)^m (3-1)^{n-m} = 2^{n-m}$, which is less than $2^n$. 

As to sharpness, see the next example.
\proofend

\begin{example} \label{curves1}
{\rm Let $\gamma_i \subset \R^2,\ i=1,\ldots,n$, be  smooth strictly convex curves. Then 
$$
L=\gamma_1\times\ldots\times\gamma_n\subset \R^2\times\ldots\times\R^2=\R^{2n}
$$
is a Lagrangian manifold. The product of the exteriors of all the curves is the domain that is covered by the tangent spaces of $L$ with multiplicity $2^n$. In terms of the local generating function, one may take
$$
F=q_1^3+q_2^3+\ldots + q_n^3.
$$
 }
\end{example}
 
\begin{example}
\label{hyp}
{\rm  Let us revisit Example \ref{Jacob}. To simplify computations, set $a=b=1$, that is, consider the generating function $F(q_1,q_2)=q_1^2q_2+q_1q_2^2$. Then the critical point equation $d (F -  x \cdot F_q+ q \cdot y)=0$ becomes the system of quadratic equations in $(q_1,q_2)$:
\begin{equation}
\label{syst}
 \left\{
   \begin{aligned}
(q_1-x_1+q_2-x_2)^2-(q_1-x_1)^2=(x_1+x_2)^2-x_1^2-y_1\\
(q_1-x_1+q_2-x_2)^2-(q_2-x_2)^2=(x_1+x_2)^2-x_2^2-y_2.\\
   \end{aligned}
\right. 
\end{equation}

\begin{figure}[hbtp]
\centering
\includegraphics[height=3.4in]{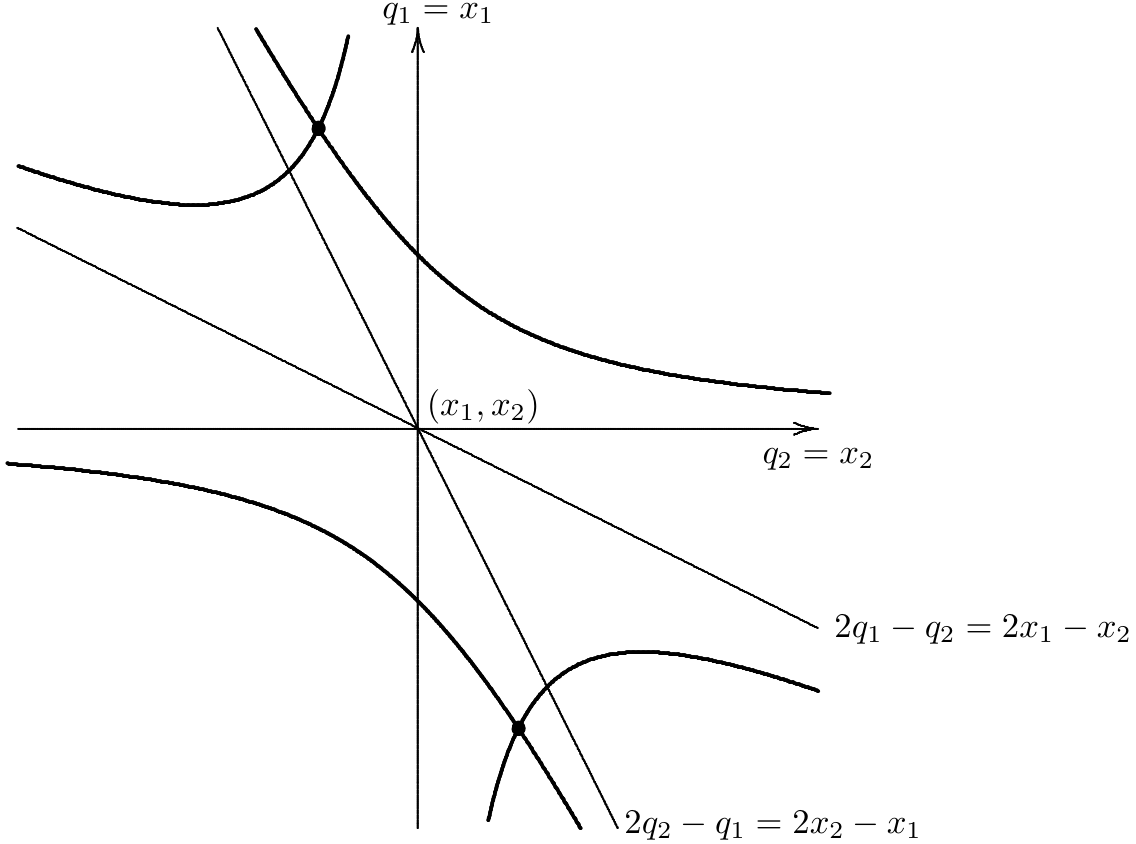}
\caption{Solution to system (\ref{syst}).}
\label{hyperbolas}
\end{figure}

The points $(x,y)\in L$ yield the trivial solution $q_1=x_1, q_2=x_2$; otherwise one has two solutions corresponding to the two intersection points of the hyperbolas with the center $(x_1,x_2)$ and alternating asymptotes in the $(q_1,q_2)$-plane, described by these equations; see Figure \ref{hyperbolas}. Thus the multiplicity is two everywhere outside of $L$.  
}
\end{example}

\section{Lagrangian outer billiards}

Consider a Lagrangian submanifold $L\subset \R^{2n}$. Two points of the ambient space are in 
the {\it outer billiard correspondence} relative to $L$ if they lie on the same tangent space to $L$ and are centrally symmetric with respect to the foot point of this tangent space.\footnote{Here we assume that this foot point is at least locally unique.} Similarly to the planar case, Theorem \ref{sympl} implies that this (possibly, multi-valued) correspondence is symplectic. 

Let us give another proof that the outer billiard correspondence a symplectic relation, by adapting the approach in \cite{Ta1,Ta2}. 

Let $\Gamma_L \subset \R^{2n} \times \R^{2n}$ be the graph of the outer billiard correspondence. Equip $\R^{2n} \times \R^{2n}$ with the symplectic structure $\omega \ominus \omega$. We want to show that $\Gamma_L$ is a Lagrangian submanifold.

Let
$$
(x,y)\quad {\rm and}\quad (\bar x,\bar y)
$$
be Darboux coordinates in the two copies of $\R^{2n}$, so that 
$$
\omega \ominus \omega = d\bar x\wedge d\bar y - dx\wedge dy.
$$
Consider the cotangent bundle $T^* \R^{2n}$ with coodinates $(q_1,q_2,p_1,p_2)$ and the standard symplectic form $\Omega=dq_1\wedge dp_1 + dq_2\wedge dp_2$.

\begin{lemma} \label{liniso}
The linear map, given by the formulas
\begin{equation} \label{iso}
q_1=\frac{x+\bar x}{2},\ q_2=\frac{y+\bar y}{2},\ p_1=\frac{\bar y-y}{2},\ p_2= \frac{x-\bar x}{2},
\end{equation}
is  (up to a constant factor) a symplectic isomorphism:
$$
\Omega = \frac{1}{2} \omega \ominus \omega.
$$
\end{lemma}

\proof
A straightforward calculation.
\proofend

Using (\ref{iso}), consider $\Gamma_L$ as a submanifold in $T^* \R^{2n}$. The next lemma implies that $\Gamma_L$ is Lagrangian.

\begin{lemma} \label{conorm}
$\Gamma_L\subset T^* \R^{2n}$ is the conormal bundle of $L\subset \R^{2n}$.
\end{lemma}

\proof
Two points $A,B \in \R^{2n}$ are in the outer billiard correspondence if and only if the midpoint of the segment $AB$ belongs to $L$ and the vector $B-A$ is tangent to $L$ at this midpoint. 

Let 
$A=(x,y), B=(\bar x,\bar y)$. Then
$((x,y), (\bar x,\bar y)) \in \R^{2n} \times \R^{2n}$ belongs to $\Gamma_L$ if and only if $q:=(q_1,q_2) \in L$ and $(-p_2,p_1) \in T_q L$. The last condition is equivalent to $(p_1,p_2)\in \nu_q (L)$ where   $\nu (L)$ is the conormal bundle.

Indeed, $(p_1,p_2)\in \nu_q (L)$ if and only if, for every test vector $(v_1,v_2)\in T_q L$, one has $p_1\cdot  v_1 + p_2\cdot v_2 =0$ or, equivalently,
$
\Omega((-p_2,p_1),(v_1,v_2))=0.
$ 
That is, $(-p_2,p_1)$ belongs to the symplectic orthogonal complement to $T_q L$. Since $L$ is Lagrangian, this orthogonal complement coincides with $T_q L$.
\proofend

\begin{remark} \label{gen}
{\rm 
One can extend the definition of the outer billiard correspondence to isotropic and to coisotropic submanifolds in linear symplectic space. Let $M\subset \R^{2n}$ be such a submanifold. Then points $A$ and $B$ are in the outer billiard correspondence with respect to $M$ if the midpoint of the segment $AB$ belongs to $M$, and the vector $B-A$ belongs to the symplectic orthogonal complement to the tangent space of $M$ at this midpoint
 (for a coisotropic manifold $M$, this implies that the line $AB$ is tangent to $M$).

The proof that this is a symplectic correspondence remains the same. The case of hypersurfaces was studied in \cite{Ta1,Ta2,Ta3}. If the hypersurface is strictly convex, this correspondence is a  symplectomorphism of its exterior.
}
\end{remark}

The following approach is a version of the ``g\' eod\' esiques bris\' ees'' techniques in symplectic topology;  it was used in \cite{Ta1,Ta2} and was borrowed from  \cite{Gi}.

Let $L\subset \R^{2n}$ be a compact Lagrangian manifold.
A cyclically ordered collection of points $z_1,z_2,\ldots,z_k$ in $\R^{2n}$ is called a $k$-periodic orbit of the outer billiard correspondence if every two consecutive points are in this correspondence relative to $L$, and there is no backtracking:
$z_{i-1} \neq z_{i+1}$ for all $i=1,\ldots,k.$
We identify periodic orbits that differ only by a cyclic permutation of the points or by orientation reversal, that is, we consider the  orbits of the dihedral group $D_k$ as single $k$-periodic orbits.

In the next theorem we assume that the closed Lagrangian manifold $L$ is non-degenerate, that is, does not contain a piece of a Lagrangian affine  space. 
Equivalently, the generating function at each point is not purely quadratic.

\begin{theorem} \label{per}
For every odd $k>1$, the outer billiard correspondence has $k$-periodic orbits.
\end{theorem}

\proof
Consider the product $({\bf R}^{2n})^{k} \times ({\bf R}^{2n})^{k}$ with the symplectic structure $\omega_1^k \ominus \omega_2^k$. Let $z_1,\ldots, z_k$ 
be  coordinates in the first, and ${\bar z_1}, \ldots, {\bar z_k}$ in the second factor 
(each $z$ and ${\bar z}$ is a vector in ${\bf R}^{2n}$). Consider two submanifolds:
$$\Gamma^k_L=\Gamma_L\times\ldots\times \Gamma_L\ (k\ {\rm times}); \quad
C = \{{\bar z_i} = z_{i+1}, i=1,\ldots,k \}.$$
The intersection $\Gamma^k_L \cap C$  consists of $k$-periodic orbits of the outer billiard relation, possibly, backtracking ones. 

As before, $({\R}^{2n})^{k} \times ({\R}^{2n})^{k}$ is symplectomorphic to $( T^* {\R}^{2n} )^k$,
and $\Gamma^k_L$ is the conormal bundle of $L^k = L \times\ldots \times L \subset ({\R}^{2n})^k$.
Let $q_1, ... , q_k$ be space and $p_1, ... , p_k$ be momenta coordinates. 

The submanifold $C$ is a linear subspace 
in $( T^* {\bf R}^{2n} )^k$, and a direct computation shows that, 
for odd $k$, it is the graph of the differential of the following 
quadratic function: 
$$\Phi (q_1,\ldots,q_k) = \sum_{1\le i < j\le k} (-1)^{i+j}\ \omega (q_i, q_j),$$ 
where $q_i\in \R^{2n}$.
Thus the points of $\Gamma^k_L\, \cap\, C$ are critical points of the restriction of $\Phi$ on $L^k$. 

Note that the function $\Phi$ is invariant under the cyclic permutations of its arguments and changes sign if the cyclic order is reversed. 

The problem is that some critical points may correspond to backtracking orbits. According to Lemma \ref{liniso}, this happens when $q_i=q_{i+1}$ for some $i$. Let us show that the maxima (or minima, which belongs to the same orbit of the dihedral group) of $\Phi$ correspond to genuine $k$-periodic orbits.

\begin{lemma} \label{degen}
In any neighborhood of an arbitrary point $q\in L$, there exist points $q_1, q_2 \in L$ such that $\omega (q_1-q, q_2-q) \neq 0$.  
\end{lemma}

\proof
Without loss of generality, assume that $q$ is the origin, and $T_q L$ is the $x$-space of the Darboux $(x,y)$-coordinate system. Let $F$ be the generating function of $L$. Working toward contradiction, assume that $\omega (q_1-q, q_2-q) = 0$ for all points $q_1, q_2 \in L$ close to the origin.

Then, in the vector notation,
$$
q_1 = (u, F_u (u)),\  q_2 = (v, F_v (v)),\ \ {\rm where}\ \ u,v \in \R^n,
$$
and 
$$
\omega (q_1, q_2) = \sum_{j=1}^n (u_j F_{v_j} (v) - v_j F_{u_j} (u)).
$$

Fix a close to the origin vector $v\in\R^n$ so that all its components are equal to zero, except $v_i\ne0$. Since  $\omega (q_1, q_2) = 0$, we have
$F_{u_i} (u) = \displaystyle\frac1{v_i}\sum u_j F_{v_j} (v)$, that is, the partial derivative $F_{u_i}$ is a linear function of $u$. Since this holds for all $i$, the function $F(u)$ is quadratic, in contradiction with the assumption that $L$ is non-degenerate.
\proofend

Assume that $(q_1,q_2,\ldots,q_k)$ is a maximum of $\Phi$ and $q_1=q_2=q$. Consider a perturbation $(\bar q_1,\bar q_2,\ldots,\bar q_k) \in L \times \ldots \times L$ with 
$$
\bar q_1 = q +u,\ \bar q_2 = q + v,\ \bar q_i = q_i\ {\rm for}\ i \ge 3.
$$
We calculate the change of the function $\Phi$:
\begin{equation} \label{delta}
\Delta \Phi  = \omega \left(v-u, q + \sum_{3\le j \le k} (-1)^j q_j\right) - \omega (u,v).
\end{equation}

Note that $\Delta \Phi$ is skew-symmetric as a function of $u$ and $v$. If one can find $u,v$ such that $\Delta \Phi \ne 0$ then, possibly after swapping $u$ and $v$, one can make  $\Delta \Phi$ positive, and $(q_1,q_2,\ldots,q_k)$ is not a maximum. Thus $\Delta \Phi =0$ for all $(u,v)$.

Note that (\ref{delta}) contains two terms, linear and quadratic in $u$ and $v$. If it vanishes identically, then both terms vanish. In particular, $\omega (u,v)=0$ for all sufficiently small $u,v$, which contradicts Lemma \ref{degen}.
\proofend

\begin{remark} {\rm In Theorem \ref{per}, one does not exclude multiple orbits (in the case when $k$ is not a prime number), so the theorem is really interesting only when $k$ is prime. Allowing multiple orbits, one can (formally) strengthen the statement: for every $k$ that is not a power of 2, there exist $k$-periodic orbits.}
\end{remark}

The actual number of periodic orbits is probably considerably greater than one.

\begin{example} \label{curves2}
{\rm Revisit Example \ref{curves1}: let $\gamma_i \subset \R^2,\ i=1,\ldots,n$, be closed smooth strictly convex curves, and 
$$
L=\gamma_1\times\ldots\times\gamma_n\subset \R^2\times\ldots\times\R^2=\R^{2n}
$$
be a Lagrangian manifold.
A set $z_1,z_2,\ldots,z_k$ in $\R^{2n}$ is a $k$-periodic outer billiard orbit if and only if, for each $i$, its projection on $i$th factor $\R^2$ is a $k$-periodic outer billiard orbit about $\gamma_i$.

By Birkhoff's theorem, for every $k\ge 3$ and every rotation number $r\leq k/2$, one has at least two such orbits. This already gives a number of $k$-periodic orbits relative to $L$ of order $k^n$. For comparison, the sum of Betti numbers of $L^k$ is $2^{kn}$.

Consider the case $k=3$. By Birkhoff's theorem, each $\gamma_i$ possesses at least two 3-periodic orbits, each can be ordered in six ways. This gives $12^n$ ordered 3-periodic outer billiard orbits relative to $L$, and factorizing by the group $S_3$, one obtains $2\cdot 12^{n-1}$ 3-periodic outer billiard orbits.
}
\end{example}

As we already mentioned, outer billiards about strictly convex hypersurfaces in symplectic spaces were studied in \cite{Ta1,Ta2,Ta3}. It is proved in \cite{Ta3} that, in $\R^{2n}$, the number of distinct 3-periodic orbits is not less than $2n$, and this estimate is sharp. 

\bigskip
{\bf Acknowledgments}. The authors are grateful to ICERM for the inspiring   and welcoming   atmosphere. The second author was supported by NSF grant DMS-1510055. Many thanks to the referee for the useful criticism.

\end{document}